\documentclass[english]{article}
\usepackage[T1]{fontenc}
\usepackage[latin9]{inputenc}
\usepackage{amsmath}
\usepackage{graphicx}
\usepackage{amssymb}
\usepackage{esint}



\usepackage{babel}

\begin{document}

\title{Introduction to generalised C\'{e}saro convergence I}

\author{Richard Stone}
\maketitle
\begin{abstract}
This is the first in a set of three papers providing an introduction to \textit{generalised} C\'{e}saro convergence. We start with traditional C\'{e}saro methods for extending classical convergence and further generalise these to allow the calculation of limits/sums for a much broader class of divergent sequences/series. These provide a constructive means of analytic continuation of functions of a complex variable and we give many examples. Future sets of papers will use these methods to derive new results (and re-derive many existing results) in areas including analytic number theory; the theory of the Riemann zeta function; reversal of order of summation; exponential sums; classical integration; Taylor series and Mellin transforms; asymptotic analysis; and a number of others.
\end{abstract}

\section{Introduction}

In this and our next two papers we present a generalisation of traditional notions of C\'{e}saro convergence. This represents an extension of existing theory regarding divergent sequences/series, and so can be understood as a set of methods and results for undertaking the analytic continuation of complex functions defined formally by such classically divergent sequences or series. While we focus on C\'{e}saro convergence, the generalisation described here extends immediately to other methods for treating divergent series, such as Borel convergence/summation.

\subsection{Future applications}

The purpose of this set of papers is to establish the key results and methods which together constitute a framework of \textit{generalised} C\'{e}saro convergence. We do this not only for its own sake, but also to provide the foundations for three subsequent sets of papers, which explore new ideas and results (as well as recasting many existing ones in new ways) in a variety of areas.

The first such set introduces C\'{e}saro arrays, which significantly expand our capacity to reverse order of summation in analysing functions. We use such arrays to derive both the asymptotics and the remarkable structure of the singularities of the function $H(z)=\frac{1}{2} + \sum_{n=1}^{\infty}\textrm{e}^{-\pi\textrm{n}^{2}\textrm{z}^{2}}$ (which famously lie at the heart of Hardy's proof that there are infinitely many non-trivial Riemann zeros on the line $Re(z)=\frac{1}{2}$); to derive many results (both old and new) regarding finite exponential sums; and to demonstrate a direct connection between such singularities on the one hand, and such algebraic results on the other.

In the second set we introduce the notion of generalised root identities, which can be viewed as generalising Hadamard's theorem. We use these, among many other applications, to derive a new family of integral identities for the argument of the Riemann zeta function, $S(T)$, conditional on the Riemann hypothesis.

In the third set, we introduce the notion of Taylor-series-to-the-left, and associated methods and results for integration, Mellin transforms, and asymptotic analysis, including local-to-global inference (or "seeing to the edge of the universe with a microscope").\footnote{These latter results extend some famous theorems of Ramanujan.}

In all these aspects - both in this last set of papers and in all the preceding sets - we show that the generalised C\'{e}saro convergence framework, rather than classical convergence, is both useful and, in many ways, actually the "right" (i.e. natural and simplifying) setting for their study.

\subsection{Apologia and a quick note concerning style}

Why have we just given a brief survey of a collection of additional, future topics and papers in section 1.1? By way of apology, the honest answer is that we think we face a fundamental, practical challenge - one which almost amounts to a "Catch-22".

Despite publishing a single paper introducing the concept of generalised C\'{e}saro convergence in the late Middle Ages ([1], 2001), the concepts and ideas expounded therein regrettably failed to inflame the imagination of the mathematical community. As such, we fear that there is likely to be limited appetite among journals for publishing the sort of refinement and extension of these ideas contained in these first three papers on their own. On the other hand, without this foundation, it would be impossible to publish the subsequent sets which contain most of the applications and results which make the foundation worth building.

As such, our hope is to break this "Catch-22" by publishing here on the Arxiv. But this is also why, in doing so, we have taken the liberty of briefly shining our torch forward, in section 1.1, onto these future payoffs - in the hope of convincing readers that there are new, beautiful and intriguing results that will flow from establishing the foundations of generalised C\'{e}saro convergence and which would reward readers for their investment of time in understanding these foundations. 

As to style, in these papers we aim to present our ideas and results as straightforwardly and simply as we can, in order to maximise their accessibility (and we emphasise at the outset that we believe that all the ideas and results presented throughout them should in fact be readily accessible to a very broad swathe of mathematicians). At the same time we aim to present them in a lively and interesting way. 

By this we mean that, while always providing sufficient rigour to prove our results in the end, we also prioritise providing the heuristic rationale for them; we often focus initially on demonstrating the utility and efficacy of the methods introduced, rather than getting bogged down in machinery and technical caveats too early; and we hope to be readable!

\subsection{Introduction to this paper}

In this paper we motivate and develop the core ideas and results introduced in [1] into a clean, simple, coherent framework. 

In sections 2.1 and 2.2 we start with basic examples and traditional notions of C\'{e}saro convergence. We explain their connection to analytic continuation in section 2.3; and in section 2.4 we reframe the discussion in terms of operators. In section 2.5 we distinguish between continuous and discrete C\'{e}saro settings and in section 2.6 we then formulate the precise definition of \textit{generalised} C\'{e}saro convergence. 

To render this practical, in section 2.7 we determine the eigenfunctions and generalised eigenfunctions of the C\'{e}saro operator in the continuous setting, and we use these in section 2.8 to translate our abstract definition of generalised C\'{e}saro convergence into a direct working version. We then reconsider the two examples of divergent series which had previously been out of reach and use this working definition to evaluate them as generalised C\'{e}saro sums, showing in the process how such generalised C\'{e}saro computations proceed in practice.

These two examples arose from consideration of the Riemann zeta function, $\zeta(s)$, at two particular $s$-values, and in section 3.1 we show, with the aid of the Euler-McLaurin sum formula, how we can in fact use generalised C\'{e}saro convergence to cleanly perform the analytic continuation of $\zeta$ from its half-plane of classical convergence to all of $\mathbb{C}\setminus\left\{1\right\}$. 

The key result we develop to facilitate this (theorem 3) is a strong C\'{e}saro-asymptotic relationship which is surprisingly simple and which is itself very useful in many other applications in future papers. 

After noting in section 3.2 some elegant related generalised C\'{e}saro limits, we further extend in section 3.3 to explain the fundamental relationship between poles in analytic continuation and points at which partial-sum functions develop pure log-divergences (representing generalised eigenfunctions of the C\'{e}saro operator with eigenvalue $1$). We use this to explain the simple pole of $\zeta$ at $s=1$ and why its residue is $1$.

The working to this point in section 3 is exclusively in the continuous C\'{e}saro setting, but in section 3.4 we also consider how it would proceed in the discrete C\'{e}saro framework. 

To this end we first calculate what the eigensequences and generalised eigensequences of the discrete C\'{e}saro operator are and how they relate to binomial coefficients. This leads to the concepts of asymptotic eigensequences and "remainder summation" (which is considered in greater depth in the next paper and very widely used in future applications). 

We are then able to replicate our analytic continuation of $\zeta$ and identification of its simple pole at $s=1$, but in this discrete C\'{e}saro setting we also find a countable collection of "anomalies" or "removable singularities". We show why these develop and how to obtain corrected evaluations at these points within the discrete C\'{e}saro framework. In so doing, we also show how the C\'{e}saro limits mentioned previously change in moving from the continuous to discrete C\'{e}saro contexts.

Finally, in section 4 we briefly consider the extension of the generalised C\'{e}saro convergence framework in two directions. 

First we cover divergence at points $z_{0}\neq{\infty}$. Secondly we extend to settings where the divergences which arise are not merely power or power-log divergences, by showing how the framework we have developed can readily be adapted to almost any form of divergence by use of a suitable measure. 

After a further intriguing aside we then conclude by briefly foreshadowing the further aspects of the generalised C\'{e}saro framework which we develop in the next paper. These include its fundamental \textit{geometric} character, its scaling and dilation invariance, its natural application to the analysis of many other special functions such as $\Gamma(s)$, and its application to the derivation of functional equations for many such functions relating their values at $s$ and $1-s$.

\section{Classical C\'{e}saro convergence and its generalisation}

Methods for assigning generalised limits (sums) to divergent sequences (series) have been studied for centuries. Many approaches exist (e.g. due to C\'{e}saro, Borel etc), each applicable to a different class of sequences/series. 

In [1] we introduced a generalisation of the traditional C\'{e}saro approach, but in a way which is also immediately applicable to the other approaches of Borel etc. 

Here we aim to give a self-contained\footnote{Going back to [1] should not be necessary, although it may occasionally be helpful and no harm will come to the reader from doing so!}, coherent presentation of this generalised C\'{e}saro methodology by beginning with the classical C\'{e}saro approach at the most elementary level, and then covering its extension in [1]. At the same time, we make things much clearer and more accessible, and introduce new methods that are better suited for further extension of the framework in subsequent papers.  
 
\subsection{Examples - successful and unsuccessful}

Let us start with four examples of traditional C\'{e}saro calculations, showing instances where such methods succeed in attaching a value to classically divergent sums and others where they fail. We will use these examples as a touchstone in developing the extension to generalised C\'{e}saro convergence methods:\\
\\
\textbf{Example (1)(a):} The series $1-1+1-1+ \ldots$\\
\\
\textbf{Example (1)(b):} The series $1+1+1+1+ \ldots$\\
\\
\textbf{Example (2)(a):} The series $1-2+3-4+ \ldots$\\
\\
\textbf{Example (2)(b):} The series $1+2+3+4+ \ldots$\\

Traditional C\'{e}saro methodology consists of trying to assign a limit to a classically divergent \textit{sequence} by repeatedly averaging it until the sequence becomes convergent. C\'{e}saro summation of divergent \textit{series} proceeds by applying this approach to the sequence of partial sums (which we shall henceforth refer to as a p-sum sequence) in the usual way.

For example in (1)(a), the p-sum sequence is $1,0,1,0,1,0, \ldots$; averaging this gives the sequence $1,\frac{1}{2},\frac{2}{3},\frac{1}{2},\frac{3}{5},\ldots$ which converges classically to $\frac{1}{2}$. We say that the C\'{e}saro sum of the original divergent series is this natural average value $\frac{1}{2}$, and write that
\begin{equation*}
\sum_{n=1}^{\infty}(-1)^{n-1} = \frac{1}{2}
\end{equation*}
in a C\'{e}saro sense.

Similarly, for (2)(a), the p-sum sequence is $1,-1,2,-2,3,-3, \ldots$. Averaging once gives the sequence $1,0,\frac{2}{3},0,\frac{3}{5},0, \ldots$. As $n\rightarrow\infty$ this approaches the sequence $\ldots ,\frac{1}{2},0,\frac{1}{2},0,\frac{1}{2},0, \ldots$ and so averaging again gives a sequence that is classically convergent to the average value of $\frac{1}{4}$. Thus the C\'{e}saro sum of $1-2+3-4+ \ldots$ is $\frac{1}{4}$ after averaging twice.

By contrast, no amount of averaging allows evaluation of the divergent series in (1)(b) and (2)(b). 

In (1)(b) the p-sum sequence $1,2,3,4, \ldots$ is linearly divergent and remains so under averaging; while in (2)(b) the p-sum sequence $1,3,6,10, \ldots$ diverges quadratically and continues to do so under averaging.

Thus traditional C\'{e}saro methods cannot handle these series. Formally, however, we can evaluate them using the C\'{e}saro calculations from (1)(a) and (2)(a) as follows. 

For (1)(b), if we let $S=1+1+1+1+ \ldots$ then $\left(1+1+1+1+\ldots\right)-\left(1-1+1-1+\ldots\right)=\left(0+2+0+2+\ldots\right)$. Hence, if we were to ignore for the moment the extra zero summands (we'll come back to this in section 2.2 below), we might surmise that $S-\frac{1}{2}=2S$ and hence that $S=-\frac{1}{2}$. 

Since $1+1+1+1+\ldots$ is the formal expression for the value of the Riemann zeta function $\zeta(s)=\sum_{n=1}^{\infty}n^{-s}$ at $s=0$ and it is well-known that $\zeta(0)=-\frac{1}{2}$, this seems to give the "right" value. 

In the same way, for (2)(b), suppose we let $\tilde{S}=1+2+3+4+ \ldots$. Then $\left(1+2+3+4+\ldots\right)-\left(1-2+3-4+\ldots\right)=\left(0+4+0+8+\ldots\right)$ and so we might argue that $\tilde{S}-\frac{1}{4}=4\tilde{S}$ and hence that $\tilde{S}=-\frac{1}{12}$, which again is the "correct" value - in this case for $\zeta(-1)$. 

These heuristic calculations suggest that, while not directly amenable to C\'{e}saro averaging, perhaps traditional C\'{e}saro methods can be adapted to also handle such series. We will need, however, to be careful.

\subsection{Pitfalls and the road to their resolution}

Not only are the calculations for (1)(b) and (2)(b) above purely formal, but even the C\'{e}saro calculations for (1)(a) and (2)(a) are subject to familiar caveats regarding series which are not absolutely convergent (let alone divergent). 

For example, it is well-known that the sum of a series which is only conditionally convergent can be changed by altering the order of the summands; and in the same vein, adding zero summands to the series in (1)(a) can also change its C\'{e}saro sum. If we take instead $1+0-1+1+0-1+1+0-1 ...$ then the p-sum sequence becomes $1,1,0,1,1,0,...$ and on averaging, this converges to the value of $\frac{2}{3}$ rather than the value of $\frac{1}{2}$ obtained above. Thus we have the following main issues:\\
\\
\textbf{Issue (i)}: Whether conditionally convergent or divergent, how can the value of a series be dependent on the order of its summands, or variable under the addition of zeros, and remain meaningful?\\
\\
\textbf{Issue (ii)}: Given that ordering and location of zeros can affect sums for non-absolutely convergent series, why was it apparently successful, in the formal calculations for series (1)(b) and (2)(b) above, to assign the same value, $S$, to $0+1+0+1+0+1+\ldots$ as to $1+1+1+\ldots$; and to assign the same value to $0+1+0+2+0+3+\ldots$ as to $1+2+3+\ldots$?\\
\\
We defer the resolution of issue (ii) to our next paper where we consider the dilation-invariance properties of generalised C\'{e}saro convergence. 

As for issue (i), however, an initial resolution consists simply of "defining the problem away" by declaring that summation of a series - whether classically or by divergent series methods - consists purely of assigning a limiting value to its sequence of p-sums. If the sequence of p-sums changes (because of rearrangement of summands or insertion of zeros), we have a different p-sum sequence and so it is reasonable that it be assigned a different limiting value. 

In the example above, the natural average value of the p-sum sequence for $1-1+1-1+\ldots$ is $\frac{1}{2}$, while the natural average value of the p-sum sequence for $1+0-1+1+0-1+\ldots$ is $\frac{2}{3}$ and it is reasonable for C\'{e}saro methods to evaluate different limiting values for these different p-sum sequences (i.e. different values for these divergent sums). 

From this perspective, the sort of anxieties noted above arise from intuitive assumptions imported from the theory of absolutely convergent series. Once we stop importing such assumptions, and consider convergence methods (including those of C\'{e}saro) as simply well-defined ways of assigning limits to classes of sequences, then such concerns go away.

Still, to understand more deeply why issue (i) is not in fact problematic, and why it does not lead to confusion in practice, we need to explore this further.

\subsection{Convergence methods as methods of analytic continuation} 

The correct way to really understand the resolution of issue (i), and indeed to understand the whole field of methods for treating divergent sequences/series, is to view such methods as being concrete approaches to the analytic continuation of functions of a complex variable. 

In short, we should view each sequence/series under review as not being examined in isolation, but rather as being a representative whose limit/sum represents the value of a complex function at some point $s\in\mathbb{C}$. 

The set of sequences/series involved may all be absolutely convergent in some region of the complex plane, but divergent in another, where divergence methods must then be invoked to evaluate the unique analytic continuation of the associated function; or in some instances the sequences/series involved may not be absolutely convergent for any $s\in\mathbb{C}$, so that a divergence method needs to be applied systematically at all points in the complex plane to even define the associated function. Either way, what is critical is not any ambiguities about the "pointwise" evaluation of the sequence/series arising at any individual point $s_{0}\in\mathbb{C}$, but that the divergence calculations involved behave in an analytic fashion as a function of $s$ for $s$ in a neighbourhood of $s_{0}$.

For instance, example series (1)(a) and (2)(a) above can be viewed as divergent series arising in the analytic continuation of the eta-function, $\eta(s):=\sum_{n=1}^{\infty}(-1)^{n-1}n^{-s}$, at the points $s=0$ and $s=-1$ lying outside the half-plane of absolute convergence $Re(s)>1$ where it is classically well-defined. The series (1)(b) and (2)(b) can be viewed in the same way for the Riemann zeta-function $\zeta(s):=\sum_{n=1}^{\infty}n^{-s}$, which is classically well-defined on the same half-plane.

Viewed this way, the natural interpretation is for the summands to arise at the integer points $n\in\mathbb{Z}_{>0}$. Thus, in the first case, $\eta(0)$ becomes $1-1+1-1+\ldots$ rather than $1+0-1+1+0-1+\ldots$ and we obtain the correct value of $\frac{1}{2}$. In fact classical C\'{e}saro methods, applied systematically as $s$ varies while consistently using this interpretation of its defining series, lead to a successful analytic continuation of $\eta$ to the whole complex plane - with one averaging being required for evaluation on the strip $0\leq{Re(s)}<1$, a second averaging being required for extension to the strip $-1\leq{Re(s)}<0$, a third for the next strip and so on. Since, for any $s_{0}\in\mathbb{C}$, the same number of averagings can be applied to evaluate $\eta(s)$ for all $s$ in a neighbourhood of $s_{0}$, this extension is analytic and thus represents the correct unique analytic continuation of $\eta(s)$ from its half-plane of absolute convergence $Re(s)>1$ to all of $\mathbb{C}$.\footnote{Here we are also using the $regularity$ of the operator $P$ to cover the situation where $Re(s_{0})$ is a non-positive integer - see next subsection.}

A nagging question, however, remains. What if we chose to be perverse and not apply the natural placement of summands, but rather systematically use the one with extra zeros after each summand $+1$ and before each summand $-1$? Would we still get the correct analytic continuation of $\eta(s)$? 

Despite the error in the evaluation of $\eta(0)$ the answer is yes! Traditional C\'{e}saro averaging would still give the correct analytic continuation of $\eta(s)$ to all $Re(s)<1$ except at the points $s\in\mathbb{Z}_{\leq{0}}$, which should be viewed as a set of "removable singularities" or "anomalies". We discuss this phenomenon (including how to correct the values at these points) in detail in section 3.4 after developing our generalised C\'{e}saro methodology. As a closing observation, however, we note that it is always good to be respectful of the \textit{natural} ordering of things and this often helps avoid such inconveniences.

Moving on to the case of $\zeta$, however, we soon see that - in line with our failures in examples (1)(b) and (2)(b) - traditional C\'{e}saro averaging is insufficient to allow any extension of $\zeta(s)$. 

Following [1], we thus now reformulate these traditional C\'{e}saro methods, which sufficed for the case of $\eta$, in terms of operators. We will then use this reformulation to extend from traditional to \textit{generalised} C\'{e}saro convergence, which will allow us to push through and handle also the analytic continuation of $\zeta$.

\subsection{Operator formulation of Cesaro convergence}

Let $\mathcal{S}$ be the space of all sequences $a=\{a_{n}\}_{n=1}^{\infty}$ and let $P_{D}:\mathcal{S} \rightarrow \mathcal{S}$ be the linear (discrete) C\'{e}saro averaging operator given by 
\begin{equation}
P_{D}[a]_{n}\colon{=}\frac{1}{n} \sum_{j=1}^{n}a_{j} \quad .
\label{P_D_Definition}\end{equation} 
Then we have the following definition of traditional C\'{e}saro convergence:\\
\\
\textbf{Definition 0:} \textit{We say that a sequence} $a\in\mathcal{S}$ \textit{has traditional C\'{e}saro limit} $L$, \textit{and we write} $C_{D}lim_{n\rightarrow\infty}a_{n}=L$, \textit{if for some positive integer} $r$ \textit{the sequence} $P_{D}^{r}[a]$ \textit{converges classically to} $L$. \\
\\
It is easy to see that $P_{D}$ is a \textit{regular} operator. By this we mean that if $a$ itself is classically convergent then so is $P_{D}[a]$ with the same limit. Hence traditional C\'{e}saro convergence is a well-defined extension of the notion of classical convergence to a broader class of sequences, since it preserves classical limits and we cannot get different limits for the same sequence by applying different powers of $P_{D}$. 

Roughly speaking, it expands the class of sequences to which we can now assign a C\'{e}saro limit to include also oscillatory sequences whose amplitude growth is bounded by some power of the index $n$. As always, the C\'{e}saro sum of a series $\sum_{i=1}^{\infty}a_{i}$ is defined to be the C\'{e}saro limit of its p-sum sequence, $\{s_{k}\}$, given by $s_{k}=\sum_{i=1}^{k}a_{i}$.

\subsection{Discrete and continuous C\'{e}saro convergence}

There is a natural continuous analogue of the discrete C\'{e}saro convergence defined above, applicable to functions rather than sequences. 

Let $\mathcal{F}$ be the space of integrable functions on finite sub-intervals of $\left[0,\infty\right)$.\footnote{In [1] we considered $\mathcal{F}$ more precisely as the space of complex-valued functions on $\left[0,\infty\right)$ given by $\mathcal{F}=\{f:\int_{0}^{x}|f(t)\left(ln(t)\right)^{m}|\,dt<\infty \; for \; all \; x\geq{1} \; and \; for \; all \; m\in\mathbb{Z}_{\geq{0}}\}$. Here, however, we embrace the liberating power of vagueness and avoid getting caught up in technical machinery in order to focus on C\'{e}saro methods and calculations as quickly and directly as possible.} We define the continuous C\'{e}saro operator $P:\mathcal{F} \rightarrow \mathcal{F}$ as the averaging operator given by 
\begin{equation}
P[f](x)=\frac{1}{x}\int_{0}^{x}f(t)dt
\label{P_Definition}\end{equation} 
for any $f\in\mathcal{F}$ and we say that a function $f$ has C\'{e}saro limit $L$ (and write $Clim_{x\rightarrow\infty}f(x)=L$) if for some positive integer $r$ the function $P^{r}[f](x)$ converges classically to $L$ as $x\rightarrow\infty$. 

Then $P$ too is a regular operator (if $\lim_{x\rightarrow\infty}f(x)=L$ then $\lim_{x\rightarrow\infty}P[f](x)=L$ also) and so continuous C\'{e}saro convergence likewise extends the notion of classical convergence of functions to a broader class of functions (roughly, oscillatory functions whose amplitude growth is bounded by a power) in a well-defined way. 

In analogy with the treatment of series in the discrete framework, we can try to make sense of a classically divergent integral, $\int_{0}^{\infty}f(x)dx$, in the continuous C\'{e}saro setting by taking the integral $F(X)\colon{=}\int_{0}^{X}f(x)dx$ and considering whether $F(X)$ has a Cesaro limit, $L$, as $X\rightarrow\infty$. If so, we say that $\int_{0}^{\infty}f(x)dx=L$ in a C\'{e}saro sense.

The discrete and continuous C\'{e}saro settings are obviously not completely distinct arenas. While we cannot apply discrete C\'{e}saro methods to problems involving functions of a continuous variable, we can use the continuous C\'{e}saro framework to handle the treatment of divergent sequences/series. 

We do this by viewing a sequence $\left\{a_{n}\right\}$ as a step-wise function, $a(x)$ on $[0,\infty)$, whose value is given by $a(x)=a_{n}$ on each interval $n\leq{x}<n+1$.\footnote{There is a technical question here about exactly how $a(x)$ should be defined at the integer points and it is usually resolved by "splitting the difference", but we omit such technicalities at this stage.} For a defining series, this is equivalent to forming a p-sum \textit{function} rather than a p-sum sequence, given by $s(k+\alpha):=\sum_{i=1}^{k}a_{i}$ where this notation (which we shall use over and over in what follows) expresses arbitrary $x\in[0,\infty)$ as $x=k+\alpha$ where $k\in\mathbb{Z}_{\geq{0}}$ and $\alpha\in[0,1)$.\footnote{Since $x$ frequently already has a meaning in many problems, we more often use $X=k+\alpha$ as our continuous C\'{e}saro limiting variable, but context will always make things clear.}

Thus, for example for series (1)(a) above, the p-sum function for $1-1+1-1+\ldots$ is $s(x)$ whose value is $0$ for $x\in[2k,2k+1)$ and $1$ for $x\in[2k+1,2k+2)$, $k\in\mathbb{Z}_{\geq{0}}$. It is easy to see that $P[s](x)$ converges classically to the "average value" of $s(x)$, namely $\frac{1}{2}$, as $x\rightarrow\infty$ and so the continuous C\'{e}saro framework assigns the same value of $\frac{1}{2}$ to this divergent series as we found before using $P_{D}$.

It turns out that the continuous C\'{e}saro framework is somewhat easier to develop and use than the discrete C\'{e}saro framework we first introduced. Thus, to avoid repetition, we will focus almost exclusively on this continuous C\'{e}saro setting in what follows - both in the generalisation of the theory of C\'{e}saro convergence in the rest of this paper and for most (but not all) of the new ideas and applications of this theory in future papers. 

In this introductory paper we will, however, make a brief excursion back into the analogous  theory for the discrete framework in section 3.4 where, as promised in section 2.2, it will be helpful in shedding light on how removable singularities can arise in the course of using generalised C\'{e}saro methods to perform constructive analytic continuation.

\subsection{From classical C\'{e}saro convergence to generalised C\'{e}saro convergence}

Let us now extend the traditional notion of C\'{e}saro convergence captured in sections 2.4 and 2.5 to define \textit{generalised} C\'{e}saro convergence. 

As discussed in [1], the key insight is simply that there is no need to restrict to pure powers of $P$ when transforming functions. Since $P$ is a regular operator, any polynomial, $q(P)$, will also be regular providing it satisfies $q(1)=1$. Thus we can generalise the traditional definition of C\'{e}saro convergence as follows:\\
\\
\textbf{Definition 1:} \textit{We say that a function} $f\in\mathcal{F}$ \textit{has generalised C\'{e}saro limit} $L$, \textit{written} $Clim_{x\rightarrow\infty}f(x)=L$, \textit{if there is a regular polynomial}, $q(P)$, \textit{such that the function} $q(P)[f](x)$ \textit{converges classically to} $L$ \textit{as} $x\rightarrow\infty$.\\

As before, regularity makes this a well-defined extension of classical convergence (and of traditional C\'{e}saro convergence, where we restricted to the case $q(P)=P^{r}$). Moreover, since a product of regular polynomials is also a regular polynomial, it is also straightforward to see that generalised C\'{e}saro convergence behaves linearly, i.e. that the C\'{e}saro limits of sums, differences and constant multiples of functions are the sums, differences and constant multiples of their C\'{e}saro limits.

What does this extended definition achieve? Well, if $q(P)$ has a regular factor $\frac{1}{1-\lambda}(P-\lambda)$ for some $\lambda\in\mathbb{\mathbb{C}}$ with $\lambda\neq{1}$ then $q(P)$ will annihilate the eigenfunction of $P$ with eigenvalue $\lambda$. Likewise, if $q(P)$ has a repeated regular factor $\frac{1}{\left(1-\lambda\right)^{N}}(P-\lambda)^{N}$ for $\lambda\neq{1}$ then $q(P)$ will annihilate the first $(N-1)$ generalised eigenfunctions of $P$ with eigenvalue $\lambda$.

Thus, our generalised notion of C\'{e}saro convergence allows us to annihilate (i.e. attach C\'{e}saro limit $0$) to eigenfunctions and generalised eigenfunctions of $P$ with eigenvalue not equal to $1$; and if two functions $f$ and $\tilde{f}$ differ by a constant multiple of any such eigenfunction or generalised eigenfunction, then they must have the same C\'{e}saro limit (providing either of them possesses a generalised C\'{e}saro limit to begin with).

As such, we can now assign a C\'{e}saro limit to any given function $f\in\mathcal{F}$ as long as we can express it as a linear combination of non-unit eigenfunctions and generalised eigenfunctions of $P$ plus a residual function, $R(x)$, which can be rendered classically convergent by applying to it a sufficiently high power, $P^{n}$ (i.e. by averaging the residual component sufficiently often in the traditional C\'{e}saro fashion).

To understand the full implications of this, of course, we need to know explicitly what the eigenfunctions and generalised eigenfunctions of P are.

\subsection{Intermission - the eigenfunctions and generalised eigenfunctions of $P$}

\textbf{Theorem 1:} \textbf{(i)} \textit{The functions} $x^{\rho}, \rho\in\mathbb{\mathbb{C}}, Re(\rho)>-1$ are all eigenfunctions of $P$ in $\mathcal{F}$ with eigenvalue $\frac{1}{\rho+1}$. Each spans a one-dimensional eigenspace of $P$.\\
\\
\textbf{(ii)} For each eigenvalue $\frac{1}{\rho+1}$, the corresponding generalised eigenfunctions of $P$ are then the functions $x^{\rho}\left(ln(x)\right)^{m}, m=1,2,3, \ldots$.\\
\\
\textbf{Proof:} It is trivial that $P[{\tilde{x}}^{\rho}](x)=\frac{1}{\rho+1}x^{\rho}$ and repeated integration by parts in turn shows that $\left(P-\frac{1}{\rho+1}\right)^{m+1}[{\tilde{x}}^{\rho}\left(ln(\tilde{x})\right)^{m}](x)=0$. QED.\\

Since we know we can assign generalised C\'{e}saro limit $0$ to eigenfunctions and generalised eigenfunctions of P, and since for $Re(\rho)\leq-1$ we already know that $x^{\rho}\left(ln(x)\right)^{m}\rightarrow0$ classically as $x\rightarrow0$, it follows that the move from traditional to generalised C\'{e}saro convergence means we can now assign limit $0$ also to all functions of the form $x^{\rho}$ or $x^{\rho}\left(ln(x)\right)^{m}, m\in\mathbb{Z}_{>1}$ for any $\rho\in\mathbb{\mathbb{C}}, \rho\neq{0}$.\\
\\
\textbf{The case of eigenvalue 1:} What about the eigenfunctions and generalised eigenfunctions with eigenvalue 1, which had to be excluded above because the regular factor $\frac{1}{1-\lambda}(P-\lambda)$ becomes singular as $\lambda\rightarrow1$? 

Well, the eigenfunction with eigenvalue 1 corresponds to $\rho=0$, i.e. the constant function $1$. Since this clearly has classical limit 1 rather than 0, it is just as well that it needs to be excluded and cannot be annihilated by any regular polynomial in $P$. 

As for the generalised eigenfunctions with eigenvalue 1, these are then the pure powers of the log-function, i.e. $\left(ln(x)\right)^{m}$. These cannot be assigned any generalised C\'{e}saro limit\footnote{Again this is just as well since $P[\ln \tilde{x}](x)=\ln x - 1$; no generalised C\'{e}saro limit, $L$, should be able to be assigned to $\ln x$ since this would then require $L=L-1$, a contradiction.} and we will see in future sections that when they arise in the course of analytic continuation, they signal the presence of poles at the complex value giving rise to the pure log-divergence.

\subsection{Generalised C\'{e}saro convergence}

Combining the observations of the last two subsections, we can then summarise what the extension from traditional to generalised C\'{e}saro convergence achieves in the following working version of theorem 1:\\
\\
\textbf{Theorem 2:} \textit{Suppose} $f\in\mathcal{F}$ \textit{can be written as} $f(x)=\sum_{j=1}^{N}c_{j}x^{{\rho}_{j}}\left(ln(x)\right)^{{m}_{j}}+R(x)$ \textit{for some collection of constants} $c_{j}\in\mathbb{C}$, $\rho_{j}\in\mathbb{C\setminus}\{0\}$ \textit{and} $m_{j}\in\mathbb{Z_{\geq\textrm{0}}}$, \textit{and some remainder function} $R(x)$ \textit{satisfying} $P^{r}[R](x)\rightarrow L$ \textit{classically as} $x\rightarrow\infty$ \textit{for some non-negative integer} $r$. \textit{Then} $Clim_{x\rightarrow\infty}f(x)=L$.\\
\\
\textbf{Notation [C\'{e}saro convergence and strong C\'{e}saro convergence]:}\\ \textbf{(i)} In addition to the notation $Clim_{x\rightarrow\infty}f(x)=L$ we shall also write $f(x)\overset{C}{\rightarrow}L$ or $f(x) \overset{C}{\sim} L$ to denote that the function $f$ converges in a generalised C\'{e}saro sense to $L$.

Functions which converge just via traditional C\'{e}saro averaging, i.e. via a pure power of $P$ without the need to remove any power or power-log divergences, will remain of interest even after this generalisation. Their divergence is in some sense better-controlled than that of a function having power or power-log divergences and it will prove critical to distinguish this sub-category of divergent functions in future working.\\
\textbf{(ii)} Thus we also introduce the  notation that $f(x) \overset{C}{\simeq} L$ and say that $f$ is \textit{strongly} C\'{e}saro convergent to $L$ if $f$ converges classically to $L$ via application of a regular polynomial $q(P)$ which is a pure power, $q(P)=P^{r}$.\\
\\
\textbf{Examples (1)(b) and (2)(b) revisited:} With theorem 2 we can now handle the two examples (1)(b) and (2)(b) - for the divergent series defining $\zeta(0)$ and $\zeta(-1)$ respectively - which were previously inaccessible to traditional C\'{e}saro methods because they had linear and quadratic power divergences. Using the notation $x=k+\alpha$ introduced earlier:\\
\\
\textbf{(i) [Example (1)(b)]:} The p-sum function for the series $\sum_{n=1}^{\infty}1$ is 
\begin{equation*}
s(x)=k=x-\alpha \quad .
\end{equation*} 
Since the saw-tooth function $R(x)=\alpha$ clearly satisfies $P[R](x)\rightarrow\frac{1}{2}$ as $x\rightarrow\infty$ it follows immediately from theorem 2 that $Clim_{x\rightarrow\infty}s(x)=-\frac{1}{2}$ via the regular polynomial $q(P)=2(P-\frac{1}{2})P$, and this agrees with the well-known value for $\zeta(0)$ as discussed previously.

The full generalised C\'{e}saro break-up of $s(x)$ in this case into three pieces - its limiting value $-\frac{1}{2}$; the eigenfunction divergence $x^{1}$; and a saw-tooth residual function with strong C\'{e}saro limit $0$ - is illustrated in the following figure:

\includegraphics[height=9cm, width=9cm]{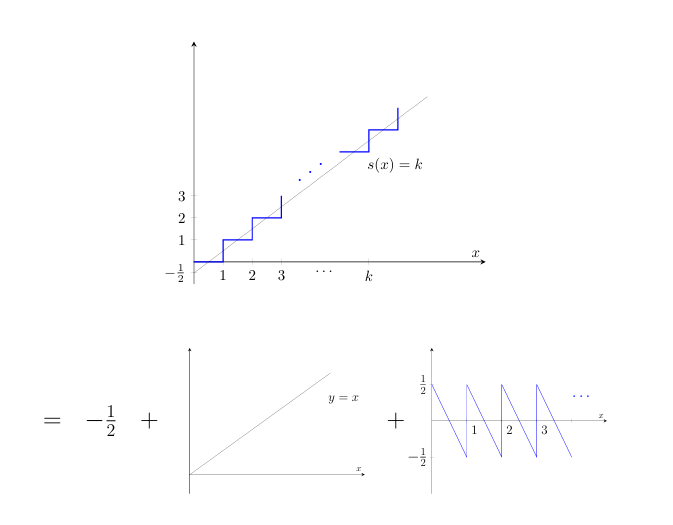}
\\
\textbf{(ii) [Example (2)(b)]:} For the series $\sum_{n=1}^{\infty}n$ formally defining $\zeta(-1)$ we have
\begin{equation}
s(x)=\frac{1}{2}(k^{2}+k)=\frac{1}{2}(k+\alpha)^{2}+(\frac{1}{2}-\alpha)k-\frac{1}{2}\alpha^{2}=\frac{1}{2}x^2+R(x)
\end{equation}
where $R(k+\alpha)=(\frac{1}{2}-\alpha)k-\frac{1}{2}\alpha^{2}$. Now \begin{eqnarray*}
P[R](k+\alpha) & = & \frac{1}{k+\alpha}\left\{\begin{array}{cc}\sum_{j=0}^{k-1}\left(\frac{1}{2} - \int_{0}^{1}\beta\,d\beta\right) \cdot j - \sum_{j=0}^{k-1}\left(\frac{1}{2}\int_{0}^{1}\beta^{2}\,d\beta\right)\\
 \\ 
 +\left(\int_{0}^{\alpha}\left(\frac{1}{2}-\beta\right)\,d\beta\right) \cdot k - \frac{1}{2}\int_{0}^{\alpha}\beta^{2}\,d\beta\end{array}\right\}\\
  \\
 & = & (-\frac{1}{6}+\frac{\alpha}{2}-\frac{\alpha^{2}}{2}) + o(1)
\end{eqnarray*}
and, since it is trivial to see that $P[\alpha^{n}](k+\alpha)\,\rightarrow\,\frac{1}{n+1}$ for any $n\in\mathbb{Z}_{\geq{0}}$, so $P^{2}[R](k+\alpha)\,\rightarrow\,-\frac{1}{12}$. Thus $Clim_{x\rightarrow\infty}s(x)=-\frac{1}{12}$ via the regular polynomial $q(P)=\frac{3}{2}(P-\frac{1}{3})P^{2}$, again in agreement with the well-known value for $\zeta(-1)$.\\

In fact, since we can now handle power and power-log divergences, generalised C\'{e}saro convergence methods actually allow us to concretely perform the analytic continuation of $\zeta(s)$ from $Re(s)>1$ to all of $\mathbb{C\setminus}\{1\}$ and to show that a simple pole with residue $1$ arises at $s=1$.

Since it is such a good "test lab" for demonstrating such analytic continuation and detection of singularities by generalised C\'{e}saro methods, we will devote section 3 to this generalised C\'{e}saro treatment of $\zeta$, noting that the same approach will apply equally well to other Dirichlet series. 

To facilitate this working, however, we conclude this section by recalling the Euler-McLaurin sum formula. This will be the critical tool for deriving from such Dirichlet series the asymptotic p-sum functions to which we can apply generalised C\'{e}saro analysis.

\subsection{The Euler-McLaurin sum formula}

We use the following working formulation (see [2]) without getting too side-tracked by detailed estimates:\\
\\
\textbf{Theorem (Euler-McLaurin sum formula)}. \textit{Suppose that} $f\in C^{\infty}\left(0,\infty\right) \cap L_{loc}^{1}\left[0,\infty\right)$ \textit{and that} $f$ \textit{and its successive derivatives form an asymptotic scale. Then we have}
\begin{equation}
\sum_{n=1}^{k}\,f(n)\,\sim\,\int_{0}^{k}\,f(x)\,dx\,+\,C_{f}\,+\,\frac{1}{2}f(k)\,+\,\sum_{r=2}^{\infty}\,\frac{B_{r}}{r!}\,f^{\left(r-1\right)}(k) \quad .
\label{EulerMcLaurin}\end{equation} \textit{Here} $C_{f}$ \textit{is a constant, the} $B_{r}$ \textit{are the alternating Bernoulli numbers} $B_{1}=-\frac{1}{2}$\textit{,} $B_{2}=\frac{1}{6}$\textit{,} $B_{3}=0$\textit{,} $B_{4}=-\frac{1}{30}$\textit{,} $B_{5}=0$\textit{,} $B_{6}=\frac{1}{42}, \ldots$ \textit{and the expression is asymptotic in the usual sense that truncating the infinite sum at any point yields a remainder which can be estimated in little-o terms by the last term retained.}\\
\\
Applied to the function $f(x;s):=\,x^{-s}$ this gives the p-sum function for the series, $\sum_{n=1}^{\infty}n^{-s}$, defining $\zeta(s)$ as
\begin{equation}
s_{\zeta,s}(k+\alpha) = \frac{k^{1-s}}{1-s}\,+\,C_{\zeta,s}\,+\,\frac{1}{2}k^{-s}+\,\sum_{r=2}^{\infty}(-1)^{r-1}\frac{B_{r}}{r!}s(s+1) \cdots (s+r-2)k^{-s-r+1}\label{psumzeta}\end{equation}
Since this is an \textit{asymptotic} series in $k$ as $k\rightarrow\infty$ it follows that for any given $s_{0}$ we can truncate the expression on the RHS after a finite number of terms, with a remainder which is $o(1)$ as $k\rightarrow\infty$ and can therefore be neglected in evaluating $Clim_{x\rightarrow\infty}s_{\zeta,s_{0}}(k+\alpha)$. By including an extra term in the RHS if required, we can ensure this $o(1)$-estimate is locally uniform in a neighbourhood of $s_{0}$.

\section{Generalised C\'{e}saro methods and the zeta function}

From equation \ref{psumzeta} we see that the series defining $s$ is absolutely convergent for $Re(s)>1$ with $\zeta(s)=C_{\zeta,s}$ on this half-plane.

\subsection{The analytic continuation of $\zeta$ via generalised C\'{e}saro methods}

Consider next the analytic continuation of $\zeta$ to the strip $0<Re(s)\leq1, s\neq{1}$. Here, equation \ref{psumzeta} tells us that
\begin{equation}
s_{\zeta,s}(k+\alpha)\,=\,\frac{k^{1-s}}{1-s}\,+\,C_{\zeta,s}\,+\,o(1)\label{psumzeta0to1}\end{equation} 
with the $o(1)$-estimate uniform on a neighbourhood of any given $s$. But for any $0\leq{Re(\gamma)}<1$ we have that $k^{\gamma}\sim{x^{\gamma}}$ classically, in the sense that $(x^{\gamma}-k^{\gamma})\rightarrow0$ classically as $x\rightarrow\infty$ (noting that $x^{\gamma}=(k+\alpha)^{\gamma}=k^{\gamma}+\gamma k^{\gamma-1}\alpha + ... = k^{\gamma} + o(1)$, again locally uniformly).

Thus for $0<Re(s)\leq1, s\neq{1}$, we have 
\begin{equation*}
s_{\zeta,s}(k+\alpha)=\frac{x^{-s+1}}{1-s}+C_{\zeta,s}+o(1)
\end{equation*} 
locally uniformly and so 
\begin{equation*}
Clim_{x\rightarrow\infty}s_{\zeta,s}(x)=C_{\zeta,s}
\end{equation*} 
via the regular polynomial $q(P;s)=\left(\frac{2-s}{1-s}\right)\left(P-\frac{1}{2-s}\right)$. 

The extension of $\zeta$ to this strip via generalised C\'{e}saro convergence thus continues to be via the same formula
\begin{equation}
\zeta(s)\,=\,C_{\zeta,s}\,=\,lim_{k\rightarrow\infty}\left(\sum_{j=1}^{k}j^{-s}-\frac{k^{1-s}}{1-s}\right)\end{equation} 
and the local uniformity of the $o(1)$-estimation above, together with the local analyticity of $q(P;s)$ in $s$ for $s\neq{1}$ ensures that this extension is the unique analytic continuation of $\zeta$ to the extended domain $Re(s)>0, s\neq{1}$.

To extend next to the strip $-1<Re(s)\leq{0}$ we need to retain an extra term in equation \ref{psumzeta}, namely
\begin{equation}
s_{\zeta,s}(k+\alpha)\,=\,\frac{k^{1-s}}{1-s}\,+\,\frac{1}{2}k^{-s}\,+\,C_{\zeta,s}\,+\,o(1).
\label{psumzeta-1to0}\end{equation}
We no longer have that $k^{\gamma}\sim{x^{\gamma}}$ classically when $1\leq{Re(\gamma)}<2$, as is the case for the leading power ($\gamma=1-s$) on this strip. But we do have that
\begin{equation}
x^{\gamma}\,\overset{C}{\simeq}\,k^{\gamma}+\frac{1}{2}\gamma k^{\gamma-1}\end{equation} in that $P\left[\tilde{x}^{\gamma}-\left(\tilde{k}^{\gamma}+\frac{1}{2}\gamma \tilde{k}^{\gamma-1}\right)\right](x)\rightarrow 0$ classically as $x\rightarrow\infty$. This result follows from the Taylor series expansion that $x^{\gamma}=\left(k+\alpha\right)^{\gamma}=k^{\gamma}+\gamma k^{\gamma-1}\alpha+o(1)$ and the fact that $P\left[\tilde{k}^{\gamma-1}\left(\alpha-\frac{1}{2}\right)\right](x)\rightarrow 0$ classically as $x\rightarrow\infty$ since $Re(\gamma-1)<1$.

Thus for $-1<Re(s)\leq{0}$, on applying this result to $k^{1-s}$ and noting that $k^{-s}\sim{x^{-s}}$ classically as before, we have, after simplification, that
\begin{equation}
s_{\zeta,s}(k+\alpha)\,\overset{C}{\simeq}\,\frac{x^{1-s}}{1-s}\,+\,C_{\zeta,s}\,+\,o(1) \quad .
\end{equation}
This is true locally uniformly (local uniformity being unaffected by application of $P$ in deducing this strongly C\'{e}saro asymptotic relationship). 

It follows immediately that we continue to have 
\begin{equation*}
Clim_{x\rightarrow\infty}s_{\zeta,s}(x)=C_{\zeta,s}\,=\,lim_{k\rightarrow\infty}\left(\sum_{j=1}^{k}j^{-s}-\frac{k^{1-s}}{1-s}-\frac{1}{2}k^{-s}\right)
\end{equation*} 
as the generalised C\'{e}saro extension of $\zeta$ to the strip $-1<Re(s)\leq{0}$, via the analytic (in $s$) regular polynomial $q(P;s)=\left(\frac{2-s}{1-s}\right)\left(P-\frac{1}{2-s}\right)P$. Local uniformity and the analyticity of $q(P;s)$ again ensure this is the unique analytic continuation of $\zeta$ to this next leftwards strip, i.e. to all of $Re(s)>-1, s\neq{1}$.

In [1] we continue in the same vein to show that generalised C\'{e}saro methodology applied to equation \ref{psumzeta} does indeed give the unique analytic continuation of $\zeta$ as $\zeta(s)=C_{\zeta,s}$ to all of $\mathbb{C\setminus}\{1\}$ by extending leftwards one strip at a time and arguing in the same fashion as above. 

The calculations involved become progressively more intricate with each new strip, but the key result derived in [1] in order to demonstrate this is that, in general, for any $Re(\gamma)\geq{0}$ we have the locally uniform strong C\'{e}saro asymptotic relationship:
\begin{equation}
x^{\gamma}=\left(k+\alpha\right)^{\gamma}\,\overset{C}{\simeq}\,k^{\gamma}+\frac{1}{2}\gamma k^{\gamma-1}+\sum_{r=2}^{\infty}(-1)^{r-1}\frac{B_{r}}{r!}\gamma\left(\gamma-1\right) \cdots \left(\gamma-r+1\right)k^{\gamma-r}\label{xtogammaexpansion}\end{equation} 
from which we derive the following fundamental result:\\
\\
\textbf{Theorem 3:} \textit{For any} $Re(s)\leq{1}, s\neq{1}$, \textit{we have locally uniformly the fundamental strong C\'{e}saro asymptotic relationship that}
\begin{equation}
s_{\zeta,s}(k+\alpha)\,=\,\sum_{j=1}^{k}j^{-s}\,\overset{C}{\simeq}\,\frac{x^{1-s}}{1-s}\,+\,C_{\zeta,s}\,+\,o(1)\label{keystrongCesaroAsymptotic}\end{equation} \textit{via the pure power} $P^{r}$ \textit{where} $r=Floor(Re(-s))+1$.\\

From this, the claimed analytic continuation of $\zeta(s)$ by generalised C\'{e}saro convergence methods as $C_{\zeta,s}$ follows at once, via the locally analytic regular polynomial $q(P;s)=\left(\frac{2-s}{1-s}\right)\left(P-\frac{1}{2-s}\right)P^{r}$, i.e.
\begin{equation}
\zeta(s)=\underset{x\rightarrow\infty}{\lim}\,q(P;s)\left[\sum_{j=1}^{\tilde{k}}j^{-s}\right]\left(x\right)=\underset{x\rightarrow\infty}{\lim}\,P^{r}\left[\sum_{j=1}^{\tilde{k}}j^{-s}\,-\,\frac{\tilde{x}^{1-s}}{1-s}\right]\left(x\right) \; .
\label{zetaclassicallimit}\end{equation}

\subsection{Remarks on theorem 3}

Theorem 3 is obviously the critical result facilitating the successful use of generalised C\'{e}saro convergence in deriving the analytic continuation of $\zeta$ from its half-plane of classical convergence $Re(s)>1$ to all of $\mathbb{C\setminus}\{1\}$. 

At least to us, however, the twin facts that equation \ref{keystrongCesaroAsymptotic} is a \textit{strong} C\'{e}saro asymptotic relationship, and that it contains only the single power $x^{1-s}$ on the RHS, are non-trivial. 

Moreover, it turns out that theorem 3 and equation \ref{keystrongCesaroAsymptotic} will play a key role in many other applications of generalised C\'{e}saro theory in future papers - for example, in many uses of "C\'{e}saro arrays"; in the analysis of the function $H(z):=\sum_{n=-\infty}^{\infty}exp\left(-\pi n^{2}z^2\right)$ which plays a central role in the analysis of the non-trivial zeros of $\zeta$; and so forth. 

As such, theorem 3 represents a foundational building block in the theory, and it is thus worth spending a moment to discuss its proof.

In [1] the working to prove theorem 3, and critical intermediate results like equation \ref{xtogammaexpansion}, is somewhat messy and involved. As a result it does not shed as much light on \textit{why} they are true as would be desirable. 

Summarising the working from [1] here in the briefest terms possible, in [1] the intermediate result \ref{xtogammaexpansion} is itself derived from a strong C\'{e}saro asymptotic formula expressing mixed products of the form $k^{\delta}\alpha^{r},\,\left(Re(\delta)\geq0, r\in\mathbb{Z_{\geq\textrm{0}}}\right)$ in terms of pure powers $k^{\delta}, k^{\delta-1}, k^{\delta-2}, ...$ (see lemma 12 in [1]). This is where the hardest of the C\'{e}saro computation is done and the key idea is to perform iterative "inversion at top order"\footnote{A la the theory for inversion of elliptic pseudodifferential operators in Analysis for those keeping score at home.}, followed by some careful combinatorics. Once the formula for $k^{\delta}\alpha^{r}$ has been obtained, we then combine it with the Taylor series for $\left(k+\alpha\right)^{\gamma}$ to obtain equation \ref{xtogammaexpansion}. The key relationship \ref{keystrongCesaroAsymptotic} then follows from inclusion of \ref{xtogammaexpansion} in the defining p-sum series \ref{psumzeta}, together with careful combinatorial simplification involving Bernoulli polynomials.

It turns out, however, that - contrary to what one might expect - this approach of expanding $\left(k+\alpha\right)^{\gamma}$ into a sum of terms of the form $k^{\delta}\alpha^{r}$ and then examining the C\'{e}saro asymptotics of these, is not the cleanest way of proving theorem 3 and getting to result \ref{keystrongCesaroAsymptotic}. Instead,\footnote{This is an insight of Donald Cartwright.} it turns out to be much simpler to consider directly the C\'{e}saro asymptotics of expressions of the form $x^{\delta}\alpha^{r},\,\left(Re(\delta)\geq0, r\in\mathbb{Z_{\geq\textrm{0}}}\right)$, and then to work backwards from these in equation \ref{psumzeta} after writing $k=x-\alpha$ and converting the RHS into a sum of such terms.

Since these results - theorem 3 and the interim results preceding it recast in this way - are not only central to future work, but also beautiful in their own right, we will prove them in this much cleaner and simpler fashion later in this set of papers introducing generalised C\'{e}saro convergence theory (specifically in the third and last paper in the set). 

However, we defer their re-framing and proof for now - other than to remark that, when we do undertake it, their proof will also introduce the important new concept of a "C\'{e}saro-adapted scale". This concept will recur in many of the more interesting pieces of future work, including in new results regarding the argument of the Riemann zeta function and regarding exponential sums.

This deferral keeps things neat and tidy here, and avoids a long digression which would derail the flow of this paper. For completeness, however, and as a guide to the intriguing nature of these core generalised C\'{e}saro results, we include here without proof the key conclusions concerning the C\'{e}saro asymptotic limits of the functions $k^{\delta}, k^{\delta}\alpha^{r}, x^{\delta}$ and $x^{\delta}\alpha^{r}$ as $x=k+\alpha \rightarrow \infty$:\\
\\
\textbf{Theorem 4:} \textit{For} $\delta\in\mathbb{C}$, $Re(\delta)\geq{0}$ \textit{and} $r\in\mathbb{Z_{\geq\textrm{0}}}$ \textit{we have}:\\
\begin{equation}
k^{\delta}\alpha^{r}\,\overset{C}{\rightarrow}\begin{cases}
0 & \quad if \; \delta\notin\mathbb{Z_{\geq\textrm{0}}}\\
\frac{(-1)^{n}}{n+r+1} & \quad if \; \delta=n\in\mathbb{Z_{\geq\textrm{0}}} \; .\end{cases}\label{CoreCesaroAsymptotic_1}\end{equation} 
\textit{Equivalently}
\begin{equation}x^{\delta}\alpha^{r}\,\overset{C}{\rightarrow}\begin{cases}
0 & \quad if \; \delta\neq{0}\\
\frac{1}{r+1} & \quad if \; \delta=0 \; .\end{cases}
\label{CoreCesaroAsymptotic_2}\end{equation}

\subsection{The singularity at $s=1$ for $\zeta$}

Returning to our analysis of $\zeta$, in section 3.1 we saw that $\zeta$ can be successfully analytically continued to all of $\mathbb{C\setminus}\{1\}$ using generalised C\'{e}saro convergence. What about the point $s=1$ in this approach?

Well, $s=1$ is in fact a simple pole of $\zeta$ with residue $1$ and this becomes apparent in the generalised C\'{e}saro approach in two ways. First, the p-sum function \textit{at} $s=1$ changes from equation \ref{psumzeta} to become instead
\begin{equation}
s_{\zeta,1}(k+\alpha)\,=\,\sum_{j=1}^{k}\frac{1}{j}\,=\,\ln k + \gamma+o(1)\,=\,\ln x + \gamma + o(1)
\label{psumzetaat1}\end{equation}
where $\gamma\approx0.577$ is the Euler-Mascheroni constant. 

The p-sum function at $s=1$ thus picks up a pure log-divergence and since $\ln x$ is a generalised eigenfunction of $P$ with eigenvalue $1$ so, as noted before, we cannot assign it a generalised C\'{e}saro limit. This signals the fact that $\zeta$ is singular at $s=1$ since, unlike the rest of $\mathbb{C}\setminus\{1\}$, it is a point where $\zeta$ cannot be assigned a value via generalised C\'{e}saro convergence.

Secondly, this acquisition of a generalised eigenfunction with eigenvalue $1$ in the p-sum function is naturally accompanied by the fact that the regular polynomial required to perform the analytic continuation across the boundary $Re(s)=1$ to the strip $0<Re(s)\leq{1}$ - namely $q(P;s)=\left(\frac{2-s}{1-s}\right)\left(P-\frac{1}{2-s}\right)$ - becomes singular as $s\rightarrow1$. It is this breakdown in the local analyticity of $q(P;s)$ that fundamentally indicates the presence of a singularity in $\zeta$ at this point.

Thus, overall, we identify where poles arise in the generalised C\'{e}saro framework either by observing where the p-sum function \textit{at} the point in question acquires a pure log-divergence, $\left(ln(x)\right)^{m}$ for some $m\in\mathbb{Z_{>\textrm{0}}}$; or else by noting where the regular polynomial needed for analytic continuation, $q(P;s)$, itself acquires a pole.

Can we tell what sort of singularity - simple pole, pole of higher order or removable singularity - occurs in such a circumstance, and also calculate either its residue (if a pole) or its true value (if a removable singularity) using our generalised C\'{e}saro framework? Yes, Wilberforce (if that is indeed the reader's name), we can - by retaining a careful focus on $q(P;s)$ while performing the necessary generalised C\'{e}saro convergence calculations.

For example, to understand the singularity at $s=1$ for $\zeta$, consider the limit as $s\rightarrow1$ of $(s-1) \cdot \zeta(s)$ expressed as the C\'{e}saro limit of $\left(s-1\right)s_{\zeta,s}(k+\alpha)$. We have
\begin{eqnarray*}
\underset{s\rightarrow1}{\lim}\left(s-1\right)\, \zeta(s)& = & \underset{s\rightarrow1}{\lim}\left(s-1\right)\, \underset{x\rightarrow\infty}{\lim}q(P;s)\left[s_{\zeta,s}\right]\left(x\right)\\ & = & -\underset{s\rightarrow1}{\lim}\, \underset{x\rightarrow\infty}{\lim}\left(\left(2-s\right)P-1\right)\left[s_{\zeta,s}\right]\left(x\right)\\
 & = & -\underset{x\rightarrow\infty}{\lim}\left(P-1\right)\left[s_{\zeta,1}\right]\left(x\right)\\
 & = & -\underset{x\rightarrow\infty}{\lim}\left(P-1\right)\left[\ln\tilde{x}+\gamma+o(1)\right]\left(x\right) = 1
\end{eqnarray*}
on swapping limits in the third line and noting that $P[\ln](x)=\ln x - 1$. 

It follows that for $s$ near $1$ we have $\zeta(s)=\frac{1}{s-1}\,+\,analytic$, and so we have derived using our generalised C\'{e}saro framework that $s=1$ is indeed a simple pole of $\zeta$ and that its residue is $1$.\\
\\
\textbf{Continuous vs discrete C\'{e}saro singularities:} Interestingly, while $s=1$ is the only place where a singularity develops in the treatment of $\zeta$ under the \textit{continuous} generalised C\'{e}saro framework, if we tackle $\zeta$ instead via the \textit{discrete} generalised C\'{e}saro framework, a further countable set of removable singularities arise at the points $s\in\mathbb{Z_{\leq\textrm{0}}}$ in addition to the pole at $s=1$. 

These turn out to be points where $q(P_{D};s)$ becomes singular, but without the associated p-sum sequences acquiring pure log-divergences. Instead, the p-sum sequences at these points acquire a \textit{constant}-sequence component corresponding to an eigensequence of $P_{D}$ with eigenvalue $1$ rather than a \textit{generalised} eigensequence with eigenvalue $1$. This explains why $q(P_{D};s)$ becomes singular at these points, but also why a limiting value can still be attached there under the discrete C\'{e}saro framework - and why these anomalous evaluations can be corrected to give the true values of $\zeta(0), \zeta(-1), \zeta(-2), \ldots$ using the discrete C\'{e}saro framework, provided we handle the polynomials $q(P_{D};s)$ carefully at these points.

Since this is worth illustrating, we now take a brief detour (as promised) into the details of the discrete C\'{e}saro frameworks Thi is worth doing in any case, since it is elegant in its own right and also provides a new way of looking at binomial coefficients and other well-known combinatorial objects.

\subsection{Understanding removable singularities using the discrete C\'{e}saro framework}

\subsubsection{Core preliminaries}

The analogue of theorem 1 in the discrete C\'{e}saro setting is as follows (but note certain clarifying remarks following the sketch proof):\\
\\
\textbf{Theorem 4:} \textbf{(i)} \textit{For} $\rho\in\mathbb{\mathbb{C}}, Re(\rho)>-1$ \textit{the eigensequence of} $P_{D}$ \textit{with eigenvalue} $\frac{1}{\rho+1}$ \textit{is given by the binomial sequence} $\lbrace\binom{n-1}{\rho}\rbrace_{n}$.\\
\\
\textbf{(ii)} \textit{The first generalised eigensequence of} $P_{D}$ \textit{with eigenvalue} $1$ \textit{is then the sequence given by} $a_{n}\,=\,1+\frac{1}{2}+\frac{1}{3}+\cdots+\frac{1}{n-1}$, \textit{and more generally, for} $\rho=m\in\mathbb{Z_{\geq\textrm{0}}}$, \textit{the first generalised eigensequence of} $P_{D}$ \textit{with eigenvalue} $\frac{1}{m+1}$ \textit{is the sequence given by} $a_{n}\,=\,\binom{n-1}{m}\left[1+\frac{1}{2}+\frac{1}{3}+\cdots+\frac{1}{n-m-1}\right]$.\\
\\
\textbf{(iii)} \textit{The first generalised eigensequence with eigenvalue} $\frac{1}{\rho+1}$ \textit{for} $\rho\in\mathbb{\mathbb{C}}, Re(\rho)>-1$ \textit{but} $\rho\notin\mathbb{Z}_{\geq{0}}$ \textit{can be expressed in analogous fashion to the integer case handled in (ii), but doing so will need the concept of "remainder summation" which we only introduce in the second of these introductory papers. We thus defer discussion of this here, as we also do for consideration of the higher-order generalised eigensequences of} $P_{D}$.\\
\\
\textbf{Remarks:} When $\rho=m\in\mathbb{Z_{\geq\textrm{0}}}$, the above claims are exact and require no refinement. For $\rho\notin\mathbb{Z_{\geq\textrm{0}}}$, however, the precise nature of the eigenvalue-relationship needs a little care. As the proof below reveals, the eigensequences either need to be understood as being indexed not from $n=1$ to $\infty$ but from $n=-\infty$ to $\infty$; or else they need to be understood in the sense of being "asymptotic eigenequences" per the definition introduced in [1], namely that\\
\\
\textbf{Definition 2:} $\lbrace a_{n} \rbrace$ \textit{is an asymptotic eigensequence of} $P_{D}$ \textit{with eigenvalue} $\lambda$ \textit{if} $\left(P_{D}-\lambda\right)\left[\lbrace a \rbrace\right]_{n}=o(1)$ \textit{classically as} $n\rightarrow\infty$.\\
\\
Further discussion of asymptotic sequences can be found in [1] (section 4.1 and remarks following proof of lemma 7) and we shall also briefly discuss some further subtleties following the proof of theorem 4.

For now, however, we wish to focus purely on proving theorem 4 and demonstrating how removable singularities arise in discrete generalised C\'{e}saro convergence theory, using $\zeta$ as a "test lab" for the investigation. For these purposes it suffices to note that for  $0\leq{Re(\rho)}<1$,  we have that 
\begin{equation}
\lbrace n^{\rho}\rbrace
\label{AsymptEigenseq_1}\end{equation} 
is an asymptotic eigensequence of $P_{D}$ with eigenvalue $\frac{1}{\rho+1}$; while for $1\leq{Re(\rho)}<2$  this extends so that
\begin{equation}
\left\lbrace n^{\rho}-\frac{\rho\left(\rho+1\right)}{2}n^{\rho-1}\right\rbrace \label{AsymptEigenseq_2}\end{equation}
is an asymptotic eigensequence of $P_{D}$ with eigenvalue $\frac{1}{\rho+1}$; and, if needed, we could continue strip-wise in this way with each new set of asymptotic eigensequences extending the previous set by inclusion of an extra term which is itself $o(1)$ on all the previous strips (so that, e.g. for $2\leq{Re(\rho)}<3$
\begin{equation}
\left\lbrace n^{\rho}-\frac{\rho\left(\rho+1\right)}{2}n^{\rho-1}+\frac{\rho\left(\rho+1\right)\left(\rho-1\right)\left(3\rho+2\right)}{24}n^{\rho-2}\right\rbrace
\label{AsymptEigenseq_3}\end{equation}
is an asymptotic eigensequence of $P_{D}$ with eigenvalue $\frac{1}{\rho+1}$ etc).\\
\\
\textbf{Proof of Theorem 4:} We shall use a useful trick and work with $P_{D}^{-1}$ rather than $P_{D}$, noting they have the same set of eigensequences and generalised eigensequences, just with reciprocal eigenvalues.

If we let $\lbrace t_{k}\rbrace:=P_{D}\left[\lbrace a \rbrace\right]_{k}$ then we have $t_{k}=\frac{1}{k}\lbrace a_{1}+\cdots+a_{k}\rbrace$ and so $a_{k}=kt_{k}-\left(k-1\right)t_{k-1}$. Thus $P_{D}^{-1}$ is given by
\begin{equation}
P_{D}^{-1}\left[\left\lbrace t\right\rbrace \right]_{k}=kt_{k}-\left(k-1\right)t_{k-1}.\label{P_D_Inverse}\end{equation}
It follows that
\begin{eqnarray*}
P_{D}^{-1}\left[\left\lbrace \binom{\tilde{k}-1}{\rho}\right\rbrace \right]_{k} & = & k\binom{k-1}{\rho}-\left(k-1\right)\binom{k-2}{\rho}\\
& = & k\binom{k-1}{\rho}-\frac{\left(k-1\right)!}{\left(\rho\right)!\left(k-2-\rho\right)!}\\
& = & \left(k-\left(k-1-\rho\right)\right)\binom{k-1}{\rho} = \left(\rho+1\right)\binom{k-1}{\rho}\\
\end{eqnarray*}
so that $P_{D}\left[\left\lbrace \binom{\tilde{k}-1}{\rho}\right\rbrace \right]_{k}=\frac{1}{\rho+1}\binom{k-1}{\rho}$ which proves (i).

Likewise, if we let $a_{k}=1+\frac{1}{2}+\cdots+\frac{1}{k-1}$, we have
\begin{eqnarray*}
P_{D}^{-1}\left[\left\lbrace a\right\rbrace \right]_{k} & = & k\left(1+\frac{1}{2}+\cdots+\frac{1}{k-1}\right)-\left(k-1\right)\left(1+\frac{1}{2}+\cdots+\frac{1}{k-2}\right)\\
& = & \left(1+\frac{1}{2}+\cdots+\frac{1}{k-2}\right)+\frac{k}{k-1} = a_{k}+1 \quad .\\
\end{eqnarray*}
Thus $\left(P_{D}^{-1}-1\right)\left[\left\lbrace a\right\rbrace \right]_{k}=1$ and therefore $\left(P_{D}^{-1}-1\right)^{2}\left[\left\lbrace a\right\rbrace \right]_{k}\equiv0$, or equivalently $\left(P_{D}-1\right)^{2}\left[\left\lbrace a\right\rbrace \right]_{k}\equiv0$, which proves the first claim in (ii).

Finally, if for $m\in\mathbb{Z_{\geq\textrm{1}}}$ we let $a_{k}=\binom{k-1}{m}\left[1+\frac{1}{2}+\cdots+\frac{1}{k-m-1}\right]$, then in the same fashion we get
\begin{eqnarray*}
P_{D}^{-1}\left[\left\lbrace a\right\rbrace \right]_{k} & = & k\binom{k-1}{m}\left[1+\cdots+\frac{1}{k-m-1}\right]\\
&   & \qquad -\left(k-1\right)\binom{k-2}{m}\left[1+\cdots+\frac{1}{k-m-2}\right]\\
& = & \left(k-\left(k-m-1\right)\right)\binom{k-1}{m}\left[1+\cdots+\frac{1}{k-m-2}\right]\\
&   & \qquad +\binom{k-1}{m}\frac{k}{k-m-1}\\
& = & \left(m+1\right)\binom{k-1}{m}\left[1+\cdots+\frac{1}{k-m-1}\right]+\binom{k-1}{m}\\
& = & \left(m+1\right)a_{k}+\binom{k-1}{m} \quad .
\end{eqnarray*}
Together with (i) this shows that $\left(P_{D}^{-1}-\left(m+1\right)\right)^{2}\left[\left\lbrace a\right\rbrace \right]_{k}\equiv0$, or equivalently that $\left(P_{D}-\frac{1}{m+1}\right)^{2}\left[\left\lbrace a\right\rbrace \right]_{k}\equiv0$, which proves the final claim in (ii) and completes the proof of theorem 4.\\
\\
\textbf{Comments: (i)} At leading order $\binom{k-1}{\rho}\sim k^{\rho}$ so that the eigensequences of (i) and the generalised eigensequences of (ii) constitute discrete analogues of the corresponding eigenfunctions ($x^{\rho}$) and generalised eigenfunctions ($x^{\rho}ln(x)$ etc) of $P$ in the continuous C\'{e}saro setting, with the same asymptotic behaviour as $k\rightarrow\infty$. While we referred readers to [1] for verification of the fact that the sequences in \ref{AsymptEigenseq_1}, \ref{AsymptEigenseq_2} and \ref{AsymptEigenseq_3} are asymptotic eigensequences of $P_{D}$, this could also be derived from a more detailed asymptotic expansion of $\binom{k-1}{\rho}$ - one which includes additional lower-order terms, together with derivation of C\'{e}saro asymptotic relationships for these terms within the discrete C\'{e}saro framework.

Whichever way one proceeds, it is easy to see from the form of the asymptotic eigensequences in \ref{AsymptEigenseq_1}-\ref{AsymptEigenseq_3} (and how they will continue to look in subsequent strips) that we have in general that $C_{D}lim_{k\rightarrow\infty}\left\lbrace k^{\rho}\right\rbrace = 0$ for $\rho\in\mathbb{C\setminus}\{\mathbb{Z_{\geq\textrm{0}}}\}$. 

The reason the cases of integer $\rho$ have to be excluded here is because, as we see in \ref{AsymptEigenseq_1}-\ref{AsymptEigenseq_3}, these asymptotic eigensequences pick up a \textit{constant} sequence component as they cross these integer $\rho$-values - that is, they pick up an eigensequence of $P_{D}$ with eigenvalue $1$, at which the associated regular polynomial $q(P_{D};\rho)$ required to obtain the generalised discrete C\'{e}saro limit therefore becomes singular. 

In fact, when $\rho=m\in\mathbb{Z_{\geq\textrm{0}}}$ it follows from the form of the exact eigensequences, $\binom{k-1}{m}$, in theorem 4(i), together with an elementary calculation, that $C_{D}lim_{k\rightarrow\infty}\left\lbrace k^{m}\right\rbrace = 1$ for all such integer $m$. We thus have the following lemma:\\
\\
\textbf{Lemma 1:} \textit{For any} $Re(\rho)\geq{0}$
\begin{equation}
C_{D}lim_{k\rightarrow\infty}\left\lbrace k^{\rho}\right\rbrace_{k}=\begin{cases}
1 & \quad if \; \rho\in\mathbb{Z_{\geq\textrm{0}}}\\
0 & \quad otherwise \; .\end{cases}
\label{k_to_the_rho_asympt_C_D_1}\end{equation}
This will be the core of what we rely on in the extension of $\zeta$ in the discrete C\'{e}saro framework.\\
\\
\textbf{(ii)} In the proof of theorem 4 we have glossed over one point. The relationship defining $P_{D}^{-1}$ in equation \ref{P_D_Inverse} involves $t_{k-1}$ and when $k=1$ this means considering $t_{0}$ even though our original definition of $\mathcal{S}$ in section 2.4 considered only sequences indexed from $1$. If we amend $\mathcal{S}$ to be sequences starting at index $0$ then we simply push the problem back one stage, since calculation of the action of $P_{D}^{-1}$ on an input sequence at index $0$ would then require us to have a value for $t_{-1}$; and so on. This is why we noted, in our remarks after stating theorem 4, that we really need to consider $P_{D}$ and $P_{D}^{-1}$ as operators applying to a space of sequences, $\mathcal{S}$, indexed from $-\infty$ to $\infty$, i.e. $\lbrace a_{n} \rbrace_{n=-\infty}^{\infty}$. 

This in turn leads to interesting issues. To begin with, how do we define $P_{D}\left[\lbrace a \rbrace\right]_{k}$ for $k\in\mathbb{Z_{\leq\textrm{0}}}$ (since $P_{D}$ and $P_{D}^{-1}$ need to map $\mathcal{S}$ back into itself and so need to be well-defined for such $k$)? And secondly, how do we make sense of our eigensequences $\binom{k-1}{\rho}$ at such indices $k\in\mathbb{Z_{\leq\textrm{0}}}$, in particular in a way which preserves the eigenvalue calculations for $P_{D}^{-1}$ in the proof?

It turns out we \textit{can} answer such questions and make sense of $P_{D}$ and $P_{D}^{-1}$ on this extended space of sequences, but doing so is somewhat non-trivial and involves two new ideas touched upon in those earlier remarks. 

The first of these is the notion of "remainder summation" which we mentioned in part (iii) of theorem 4. We will only get to this in the sequel to this paper where it will then allow us to make sense of the definition of $P_{D}$ for indices $k\in\mathbb{Z_{\leq\textrm{0}}}$. 

The second is the interpretation of sequences $\lbrace a_{n} \rbrace$ in $\mathcal{S}$ as being, in a sense, embedded in the complex plane. This allows us to think of $n$ adopting non-integer values and thereby to get at the definition of objects like $\binom{k-1}{m}$ when $m\in\mathbb{Z}$ by considering $\binom{s-1}{m}$ for $s\in\mathbb{C}$ and taking the limit as $s\rightarrow k$. 

Both of these are important ideas and recur extensively in future work, but we avoid any further discussion of them here, and point the reader\footnote{Yes indeed, I am addressing you young Wilberforce my boy!} to these future papers where he can see them properly defined and utilised.

\subsubsection{Removable singularities}

Let us now return to our original aim in this subsection - to see how removable singularities arise when trying to analytically continue $\zeta$ in the discrete C\'{e}saro framework.

Let $\zeta_{D}^{ext}(s)$ be the function obtained by applying the discrete C\'{e}saro framework to the defining series for $\zeta$, $\sum_{n=1}^{\infty}n^{-s}$, in the same fashion as was done in section 3.1 using the continuous C\'{e}saro framework.

Using the asymptotic series \ref{psumzeta} for the p-sum sequence of this defining series, we see that $\zeta_{D}^{ext}(s)=C_{\zeta,s}=\zeta(s)$ trivially for $Re(s)>1$. 

Likewise, from equation \ref{psumzeta0to1} in the strip $0<Re(s)\leq{1}$ it follows immediately, in light of equation \ref{k_to_the_rho_asympt_C_D_1}, that we continue to have $\zeta_{D}^{ext}(s)=C_{\zeta,s}=\zeta(s)$ on this strip, at least for $s\neq{1}$ (where the presence of a simple pole of residue $1$ can be deduced in the discrete C\'{e}saro framework in essentially identical fashion to that used in subsection 3.3 in the continuous C\'{e}saro setting).

At $s=0$, however, we immediately appear to have a problem. There we have p-sum sequence $s_{k}=k$ and so, in light of equation \ref{k_to_the_rho_asympt_C_D_1} we have $\zeta_{D}^{ext}(0)=1$, which does not agree with the correct value of $\zeta(0)=-\frac{1}{2}$.

In the same way, for the next strip $-1<Re(s)\leq{0}$, it follows immediately from the p-sum sequence \ref{psumzeta-1to0} and equation \ref{k_to_the_rho_asympt_C_D_1} that we have $\zeta_{D}^{ext}(s)=C_{\zeta,s}=\zeta(s)$ for $s\neq{0}$; but when we get to the point $s=-1$ we have p-sum sequence $s_{k}=\frac{1}{2}k^{2}+\frac{1}{2}k$, and so from equation \ref{k_to_the_rho_asympt_C_D_1} we again deduce that $\zeta_{D}^{ext}(-1)=1$ which does not agree with the correct value of $\zeta(-1)=-\frac{1}{12}$.

In fact, proceeding leftwards in the same way, one strip at a time, we find that extension of $\zeta$ from $Re(s)>1$ to all of $\mathbb{C\setminus}\{1\}$ using discrete C\'{e}saro methods gives a function $\zeta_{D}^{ext}$ which equals $\zeta(s)$ for all $s\in\mathbb{C\setminus}\{\mathbb{Z_{\leq\textrm{0}}}\}$, but which has the erroneous value $1$ at all non-positive-integer points $s$ (see [1], lemma 14 for proof of this last assertion).

These anomalous values at $s\in\mathbb{C\setminus}\{\mathbb{Z_{\leq\textrm{0}}}\}$  represent a series of removable singularities of $\zeta_{D}^{ext}$ where, as hinted, the regular polynomial $q(P_{D};s)$ used in the evaluation of $\zeta_{D}^{ext}(s)$ for $Re(s)\leq{1}$ becomes singular, but without the p-sum sequences at these points acquiring the pure log-divergences that would signal the presence of poles. 

For example, consider the first such anomalous point, $s=0$. Any punctured neighbourhood of this point extends into the strip $-1<Re(s)\leq{0}$ and so we need to be able to handle the p-sum sequences with both $k^{1-s}$ and $k^{-s}$ terms as given by equation \ref{psumzeta-1to0} in this strip (noting equation \ref{psumzeta-1to0} also applies for $Re(s)>0$ where the additional $k^{-s}$ term is classically $o(1)$). 

Rewriting this in terms of the asymptotic eigensequences of $P_{D}$ given in equations \ref{AsymptEigenseq_1} and \ref{AsymptEigenseq_2}, we have
\begin{equation}
\left(s_{\zeta,s}\right)_{k}=\frac{1}{1-s}\left(k^{1-s}-\frac{\left(1-s\right)\left(2-s\right)}{2}k^{-s}\right)+\frac{3-s}{2}k^{-s}+C_{\zeta,s}+o(1) \label{zeta_ext_D_0_psum}\end{equation}
and it follows that the regular polynomial we need to use to obtain discrete C\'{e}saro convergence to $\zeta_{D}^{ext}$ throughout such a neighbourhood is
\begin{equation}
q(P_{D};s)=-\frac{2-s}{s}\left(P_{D}-\frac{1}{2-s}\right)\left(P_{D}-\frac{1}{1-s}\right) \quad . \label{zeta_ext_D_0_regpoly}\end{equation}
Here the first factor annihilates the eigensequence $\lbrace k^{1-s}-\frac{\left(1-s\right)\left(2-s\right)}{2}k^{-s} \rbrace_{k=1}^{\infty}$ and the second the other eigensequence $\lbrace k^{-s} \rbrace_{k=1}^{\infty}$ in expression \ref{zeta_ext_D_0_psum}.

But it is now clear why the calculation of $\zeta_{D}^{ext}(0)$ was anomalous - the analyticity of $q(P_{D};s)$ breaks down at $s=0$ due to the development of a simple pole in the regularising factor $\frac{2-s}{s}$. 

However, as noted, this breakdown arises only because the p-sum sequence at $s=0$ in expression \ref{zeta_ext_D_0_psum} acquires a \textit{constant} eigensequence of $P_{D}$ with eigenvalue $1$, rather than acquiring a pure log-divergent \textit{generalised} eigensequence of $P_{D}$ with eigenvalue $1$. As such the anomalous value of $\zeta_{D}^{ext}$ at $0$ is a removable singularity rather than a pole of $\zeta_{D}^{ext}$, and as before, we can carefully allow for it using a L'Hopital's law calculation and still deduce the correct value for $\zeta(0)$ within our discrete C\'{e}saro scheme, as follows:
\begin{eqnarray*}
\zeta(0) & = & \underset{s\rightarrow0}{\lim}\;\zeta_{D}^{ext}\left(s\right)\;=\; \underset{s\rightarrow0}{\lim}\;\underset{k\rightarrow\infty}{\lim}\;q\left(P_{D};s\right)\left[\lbrace \left(s_{\zeta,s}\right)_{\tilde{k}} \rbrace\right]_{k}\\
 \\
& = & -\underset{k\rightarrow\infty}{\lim}\;\underset{s\rightarrow0}{\lim}\;\frac{\left(\left(2-s\right)P_{D}-1\right)\left(P_{D}-\frac{1}{1-s}\right)\left[\lbrace \left(s_{\zeta,s}\right)_{\tilde{k}} \rbrace\right]_{k}}{s}\\
 \\
& = & -\underset{k\rightarrow\infty}{\lim}\;\underset{s\rightarrow0}{\lim}\left\{\begin{array}{cc} \left((2-s)P_{D}-1\right)\left(P_{D}-\frac{1}{1-s}\right)\left[\left\{ \frac{d}{ds}\left(s_{\zeta,s}\right)_{\tilde{k}} \right\}\right]_{k}\\
 \\
 + \left(\frac{-1}{\left(1-s\right)^{2}}\right)\;\left((2-s)P_{D}-1\right)\left[\left\{ \left(s_{\zeta,s}\right)_{\tilde{k}} \right\}\right]_{k}\\
 \\
 - P_{D}\;\left(P_{D}-\frac{1}{1-s}\right)\left[\left\{ \left(s_{\zeta,s}\right)_{\tilde{k}} \right\}\right]_{k} \end{array}\right\}\\
 \\
& = & -\underset{k\rightarrow\infty}{\lim}\;\left\{\begin{array}{cc} \left(2P_{D}-1\right)\left(P_{D}-1\right)\left[\left\{\begin{array}{cc} -\tilde{k}\ln\tilde{k}+\tilde{k} \\
 +\left(\frac{d}{ds}C_{\zeta,s}\right)_{s=0} -\frac{1}{2}\ln\tilde{k} \end{array}\right\}\right]_{k}\\
 \\
 - \left(2P_{D}-1\right)\;\left[\left\{ \tilde{k} \right\}\right]_{k}\;-\;P_{D}\;\left(P_{D}-1\right)\left[\left\{ \tilde{k} \right\}\right]_{k} \end{array}\right\}\\
& = & -\underset{k\rightarrow\infty}{\lim}\;\left\{ \left(\frac{5}{4}-\frac{1}{4}k+\frac{1}{2}\right)-1+\frac{1}{4}\left(k-1\right) \right\} \;=\;-\frac{1}{2}\\
\end{eqnarray*}
where in the final steps here we have used that $\left(2P_{D}-1\right)\left[\left\{ \tilde{k}\ln\tilde{k} \right\}\right]_{k}=\ln k - \frac{1}{2}k + o(1)$ and $\left(P_{D}-1\right)\left[\left\{ \ln\tilde{k} \right\}\right]_{k}=-1 + o(1)$.

\section{Further extensions - Divergence approaching a point $z_{0}\neq{\infty}$; The general form of generalised convergence schemes; Where next}

In this first paper we have introduced the basic theory of generalised C\'{e}saro convergence in both its continuous and discrete frameworks; shown how it can be used in the case of $\zeta$ to perform analytic continuation from the region where the defining series classically converges to the rest of the complex plane where it becomes divergent; explored how we can detect the presence of either poles or anomalies/removable singularities within the C\'{e}saro framework; and then shown how we can calculate either their associated residues or their corrected values using generalised C\'{e}saro calculations.

We now quickly cover a few further natural extensions of this generalised C\'{e}saro theory, before briefly outlining the topics covered in the second paper in this introductory set.

\subsection{Divergence approaching a point $z_{0}\neq{\infty}$}

So far, we have focused exclusively on divergence as $X\rightarrow\infty$. In many instances, however, divergence issues arise at a finite point $z_{0}$. 

In such cases, the whole generalised continuous C\'{e}saro machinery still applies providing we simply change variables to $z_{0}+\frac{1}{X}$ and then consider $Clim_{X\rightarrow\infty}$. Under the change of variable $\epsilon = \frac{1}{X}$, power and power-log divergences in $X$($X^{\rho}$ or $X^{\rho}\left(\ln X\right)^{m}$) remain power and power-log divergences in $\epsilon$ ($\epsilon^{-\rho}$ or $\epsilon^{-\rho}\left(\ln\epsilon\right)^{m}$), so this amounts to understanding that a generalised C\'{e}saro approach allows us to "throw away" power and power-log divergences wherever they arise, using regular polynomials in a suitable C\'{e}saro operator ($P_{X}$ or $P_{\epsilon}$) associated to the relevant limiting variable. We just need to pay attention to doing so in an analytic fashion (i.e. with the regular polynomial involved being analytic in any underlying complex variable).

Thus, for instance, if $f$ is continuous on $(0,\infty)$ but potentially non-integrable at $0$ as well as $\infty$, we can nonetheless consider a generalised C\'{e}saro definition of the definite integral $\int_{0}^{\infty}f(x)\,dx$ as
\begin{equation}
\int_{0}^{\infty}f(x)\,dx:=\underset{X\rightarrow\infty}{Clim}\int_{1/X}^{X}f(x)\,dx \quad .
\label{definite_integral_Cesaro_1}\end{equation}
Technically here we should consider
\begin{equation}
\int_{0}^{\infty}f(x)\,dx:=\underset{X\rightarrow\infty,Y\rightarrow\infty}{Clim}\int_{1/Y}^{X}f(x)\,dx \label{definite_integral_Cesaro_2}\end{equation}
via separate regular polynomials in C\'{e}saro operators $P_{X}$ and $P_{Y}$,  but we will frequently be loose and simplify our calculations by taking $Y=X$. As long as we are careful not to illegally cancel generalised eigenfunctions with eigenvalue $1$ arising near $\infty$ (i.e. pure log-divergences in $X$) with generalised eigenfunctions with eigenvalue $1$ arising near $0$ (i.e. pure log-divergences in $Y$) when we do this, then we will not end up introducing errors into such C\'{e}saro calculations and no harm will come to man or beast\footnote{we rely on Wilberforce's vigilance in policing this in future papers!}.

As an example, consider the definite integral $\int_{0}^{\infty}\frac{x^{s-1}}{x+1}\,\textrm{d}x$ defining the Mellin transform of $\frac{1}{1+x}$ as a function of the complex variable $s$. Classically, the integrand diverges too rapidly to be integrable at $0$ when $Re(s)\leq{0}$ and is non-integrable at $\infty$ when $Re(s)\geq{1}$. For $0<Re(s)<1$ we obtain a well-defined holomorphic function, but analytically continuing it outside this strip in both directions requires additional work. 

Viewing $\int_{0}^{\infty}\frac{x^{s-1}}{x+1}\,\textrm{d}x$ as a generalised C\'{e}saro integral as above, however, accomplishes this analytic continuation immediately. We at once get a well-defined function of $s$ on the whole complex plane, and indeed it is straightforward to see within this continuous C\'{e}saro framework that this Mellin transform has simple poles at $s\in\mathbb{Z}_{\leq{0}}$ (where we pick up a log-divergence from integration near $0$ based on the Taylor series $\frac{1}{1+x}=1-x+x^{2}-\ldots$) and at $s\in\mathbb{Z}_{>0}$ (where we pick up log-divergences from integration near $\infty$ after asymptotically expanding $\frac{1}{x+1}$ as $\frac{1}{x}\left\{ 1-\frac{1}{x}+\frac{1}{x^{2}}-\cdots \right\}$.

This example in fact serves as a first illustration of what we believe is a general precept - one which we will  examine in detail in the third set of papers following this introductory set - namely that \textit{the continuous C\'{e}saro framework is the correct perspective from which to view the entire theory of the Mellin transform}. Viewed from this perspective the whole theory makes sense naturally and many of the usual caveats and restrictions regarding domains of convergence and regions of invertibility disappear.

We conclude by noting that if $f$ becomes singular also at a point $z_{0}\in\mathbb{R}_{>0}$, in addition to possible non-integrability at $0$ and $\infty$, then we can still define the generalised C\'{e}saro integral $\int_{0}^{\infty}f(x)\,dx$ as
\begin{equation}
\int_{0}^{\infty}f(x)\,\textrm{d}x:=\underset{X\rightarrow\infty}{Clim}\left\{ \int_{1/X}^{z_{0}-\frac{1}{X}}f(x)\,\textrm{d}x + \int_{z_{0}+\frac{1}{X}}^{X}f(x)\,\textrm{d}x \right\} \quad .
\label{definite_integral_Cesaro_3}\end{equation}
Here again we need to be careful to avoid impermissible cancellation of pure log-divergences if they arise at different points among $0, z_{0}$ and $\infty$; and we need to consider the integral more rigorously as
\begin{equation}
\int_{0}^{\infty}f(x)\,\textrm{d}x:=\underset{X_{1},\cdots,X_{4}\rightarrow\infty}{Clim}\left\{ \int_{1/X_{1}}^{z_{0}-\frac{1}{X_{2}}}f(x)\,\textrm{d}x + \int_{z_{0}+\frac{1}{X_{3}}}^{X_{4}}f(x)\,\textrm{d}x \right\} \label{definite_integral_Cesaro_4}\end{equation}
if doing so is required to avoid any such ambiguity.

Further extension of such C\'{e}saro integral definitions to encompass the case of multiple singularities in the integrand on $(0,\infty)$, or to other contours in the complex plane then follow in the same way.

\subsection{The general form of generalised convergence schemes}

In section 7 of [1] we discussed at some length how the ideas introduced here to generalise traditional C\'{e}saro convergence apply equally well to other well-known traditional convergence methods, such as Borel's scheme for handling series with exponential (rather than power) divergences. We summarise here.

In general, if $\mu=\mu(t)\,dt$ is a measure on $[0,\infty)$ such that $\mu(t)\geq{0}$ for all $t\geq{0}$ and $\int_{0}^{X}\mu(t)\,\textrm{d}t$ increases unboundedly, then we can mimic our C\'{e}saro approach and define a generalised convergence scheme adapted to $\mu$. 

We define the operator $P_{\mu}$ by 
\begin{equation*}
P_{\mu}\left[g\right](X):=\frac{1}{F_{\mu}(X)}\int_{0}^{X}g(x)\mu(x)\,\textrm{d}x
\end{equation*} 
where $F_{\mu}(X)=\int_{0}^{X}\mu(x)\,\textrm{d}x$. 

The eigenfunctions of $P_{\mu}$ are then the functions $\left(F_{\mu}(x)\right)^{\rho}$, with eigenvalue $\frac{1}{\rho+1}$, and the corresponding generalised eigenfunctions with eigenvalue $\frac{1}{\rho+1}$ are the functions $\left(F_{\mu}(x)\right)^{\rho}\left(\ln\left(F_{\mu}(x)\right)\right)^{m}$, $m\in\mathbb{Z}_{\geq{1}}$. 

In exact analogy with our C\'{e}saro approach, we say that a classically divergent function, $f$, converges to limit $L$ in a generalised $P_{\mu}$-sense if we can find a regular polynomial, $q\left(P_{\mu}\right)$, such that $q\left(P_{\mu}\right)\left[f\right](X)\rightarrow L$ classically as $X\rightarrow\infty$. 

Hence $f$ converges to $L$ in a generalised $P_{\mu}$-sense if we can express $f(x)$ as a linear combination of eigenfunctions and generalised eigenfunctions of $P_{\mu}$ with eigenvalues not equal to $1$, plus a residual function which becomes classically-convergent to $L$ on being $P_{\mu}$-averaged sufficiently often (i.e. on application of $P_{\mu}^{N}$ for some $N\in\mathbb{Z}_{\geq{0}}$). 

The case of generalised C\'{e}saro convergence which we focus on in this and subsequent papers simply corresponds to the case of $\mu(t)\equiv1$.

An idea canvassed in [1], section 7, is that if we are considering a divergent series representing the analytic continuation of an underlying complex function, we might be able to choose $\mu(t)$ suitably so that the above $P_{\mu}$-eigenfunctions and generalised eigenfunctions match the nature of the divergences which arise, and thereby endeavour to perform the analytic continuation using the generalised $P_{\mu}$-convergence scheme. For example, where we encounter exponential divergences outside a region of classical convergence, we might attempt to perform analytic continuation using a generalised Borel-type convergence scheme based, roughly speaking, on taking $\mu(t)=e^{t}$.

We do not, however, explore this in any further detail here, since we wish to retain a clear focus purely on generalised C\'{e}saro convergence. While there is much of potential interest to examine along these lines (and these virgin lands are certainly not currently teeming with explorers), we refer the reader instead to the discussion in [1]. 

The only remaining thing we note in this regard is that there is an additional interesting discussion in [1], subsection 8.3, of the one-parameter family of such schemes, $P_{r}$, $r>-1$ given by taking $P_{r}:=P_{\mu_{r}(t)}$ where $\mu_{r}(t)=t^{r}$; in particular, of how the $P_{r}$-eigenfunctions and generalised eigenfunctions vary as $r$ varies and what happens in the limiting cases $r=\infty$ and $r=-1$.

\subsection{An intriguing aside}

We conclude the detailed content of this paper with an intriguing side-observation regarding the zeta function - one which points the way towards one of the main new ideas covered in the next paper in this introductory set on generalised C\'{e}saro convergence.

For $s\in\mathbb{Z}_{\leq{0}}$ we have that
\begin{equation}
\zeta(s)=\int_{-1}^{0}\,\sum_{j=1}^{k}j^{-s}\,\textrm{d}k \quad . \label{zeta_definite_integral_1}\end{equation}
This follows because for $s\in\mathbb{Z}_{\leq{0}}$ the sum $\sum_{j=1}^{k}j^{-s}$ is a polynomial in $k$ and theorem 4, equation \ref{CoreCesaroAsymptotic_1} says that for integer powers,
\begin{equation}
\underset{k\rightarrow\infty}{Clim}\,k^{n}=\frac{\left(-1\right)^{n}}{n+1}=\int_{-1}^{0}k^{n}\,dk \quad .
\end{equation}

Apart from its intrinsic interest, the intriguing aspect of equation \ref{zeta_definite_integral_1} is that the integral in it recasts the discrete summation variable $k$ as a continuous integration variable on $[-1,0]$. 

In our next paper, however, we will make this completely natural by introducing the notion of remainder summation in the generalised C\'{e}saro framework. This will allow us to recast the classical notion of a finite sum, $\sum_{j=1}^{k}$, instead as a function of an arbitrary complex variable, $k$. This can then be integrated as above, or differentiated w.r.t. $k$ etc. One corollary will be a generalisation of result \ref{zeta_definite_integral_1} to a corresponding result for arbitrary $s\in\mathbb{C}$.

\subsection{Where next}

Along with the concept of remainder summation, we extend the theory of generalised C\'{e}saro convergence in a number of directions in the next paper. What delights await not only Wilberforce, but the swelling army of readers now accompanying him, in the next instalment introducing generalised C\'{e}saro convergence theory? 

To begin with, as part of the extension to remainder summation, we will see that \textit{geometry} becomes crucial in C\'{e}saro analysis in very interesting ways; and related to this, we will demonstrate the dilation-invariance, and also scaling-invariance, of C\'{e}saro convergence. 

As a demonstration of the value of these extensions, we show that combining them allows us to express the classical Gamma function, $\Gamma(s)$, in a new and more natural way - and in turn to deduce certain of its more advanced and mysterious properties, such as its multiplication identities, in an almost trivial way from geometric considerations. Using the same geometric and dilation properties, we also answer a number of small questions posed in this paper.

By considering remainder summation in both directions (positive and negative) we also derive functional equations for a number of well-known functions and see that generalised geometric C\'{e}saro methods shed light on why it is that so many such functional equations share the same core structure - relating function values at $s$ and $1-s$.

\section{Acknowledgements}

We thank Professor T. Abby for his help in preparing this paper, and Christiana Stone for her help in preparing the figure.


\begin{thebibliography}{4}
\bibitem{key-1}R. Stone, \textit{Operators and Divergent Series},
Pacific Journal of Mathematics, Vol. 217, No. 2, 2004

\bibitem{key-2}G.H. Hardy, \textit{Divergent Series}, Oxford, at the Clarendon Press, 1949, MR 0030620 (11,25a), Zbl 0032.05801
\end{thebibliography}
\end{document}